\documentclass{article}

\usepackage{lipsum}
\usepackage{graphicx}
\usepackage{epstopdf}
\usepackage{algorithm,algorithmic}
\usepackage{amsopn,amsfonts,amsmath,amssymb,amsthm}
\usepackage{hyperref}

\def\xt#1{{\textbf{#1}}^{\top}}
\def\x#1{\textbf{#1}}
\def\bGamma{\boldsymbol\Gamma}

\DeclareMathOperator{\diag}{diag}
\DeclareMathOperator{\tr}{tr}
\DeclareMathOperator{\stril}{stril}

\title{Joint Approximate Diagonalization under Orthogonality Constraints}

\author{Ronald de Vlaming\footnote{School of Business and Economics, Vrije Universiteit Amsterdam, Amsterdam, The Netherlands 
  (\url{mailto:r.devlaming@vu.nl}, \url{https://research.vu.nl/en/persons/ronald-de-vlaming}).} \, and \,Eric A.W. Slob\footnote{MRC Biostatistics Unit, School of Clinical Medicine, University of Cambridge, Cambridge, United Kingdom 
  (\url{mailto:eric.slob@mrc-bsu.cam.ac.uk}).}}

\begin{document}

\maketitle

\section*{Abstract}

Joint diagonalization of a set of positive (semi)-definite matrices has a wide range of analytical applications, such as estimation of common principal components, estimation of multiple variance components, and blind signal separation. However, when the eigenvectors of involved matrices are not the same, joint diagonalization is a computationally challenging problem. To the best of our knowledge, currently existing methods require at least $O(KN^3)$ time per iteration, when $K$ different $N \times N$ matrices are considered. We reformulate this optimization problem by applying orthogonality constraints and dimensionality reduction techniques. In doing so, we reduce the computational complexity for joint diagonalization to $O(N^3)$ per quasi-Newton iteration. This approach we refer to as JADOC: Joint Approximate Diagonalization under Orthogonality Constraints. We compare our algorithm to two important existing methods and show JADOC has superior runtime while yielding a highly similar degree of diagonalization. The JADOC algorithm is implemented as open-source Python code, available at \url{https://github.com/devlaming/jadoc}. 

\noindent \textbf{Key words.} diagonalization, blind signal separation, variance component estimation

\noindent \textbf{AMS subject classifications.} 65F25, 90C53

\section{Introduction}
Given a set of $N \times N$ symmetric matrices, $\x{C}_k$ for $k=1,\ldots,K$, that are (at least) positive semidefinite (PSD) the objective of finding an $N \times N$ matrix $\x{B}$ such that $\x{B}\x{C}_k\xt{B}$ is `as diagonal as possible' for $k=1,\ldots,K$ arises naturally in estimation of common principal components \cite{pham2001}.

Assuming matrices $n_k \x{C}_k$, for $k=1,\ldots,K$, follow Wishart distributions with $n_k$ degrees of freedom and covariance matrix $\bGamma_k$ respectively, such that (\emph{i}) $\bGamma_k$ have the same eigenvectors for $k=1,\ldots,K$ and (\emph{ii}) $n_k \x{C}_k$ and $n_l \x{C}_l$ are independent for $k\neq l$, one can apply maximum likelihood estimation (MLE) to obtain estimates of these common eigenvectors. The resulting optimization problem can be rewritten as finding a matrix $\x{B}$ that minimizes the following criterion \cite{pham2001}:

\begin{equation}\label{eq:crit_pham}
\mathcal{L}\left( \x{B} \right)  = \frac{1}{2} \sum_{k=1}^K n_k \left(\log \det \diag \left( \x{B} \x{C}_k \xt{B} \right) - \log \det \left( \x{B} \x{C}_k \xt{B} \right) \right).
\end{equation}
This criterion provides a reasonable overall measure of deviation of $\x{B} \x{C}_k \xt{B}$ from diagonality for $k=1,\ldots,K$. Thus, this MLE-based criterion, as proposed by Pham \cite{pham2001}, is a natural starting point for finding a matrix $\x{B}$ that aims to diagonalize $\x{B}\x{C}_k\xt{B}$ for $k=1,\ldots,K$.

Subsequent research on this general approach to find $\x{B}$ has been fairly limited, however, possibly owing to the fact that this optimization problem is computationally challenging and has poor scalability \cite{mesloub2014}. Existing approach typically use so-called sweeps \cite{pham2001}. Within a sweep, an optimization is performed for all pairwise combinations of rows of $\x{B}$.

Importantly, relatively recently, a quasi-Newton method has been proposed to find a matrix $\x{B}$ that minimizes Equation \ref{eq:crit_pham} more efficiently \cite{ablin2018}. The great advantage of that particular approach is that all elements of $\x{B}$ are updated jointly in each iteration. However, that approach requires at least $K$ matrix multiplications of $N \times N$ matrices in order to calculate the gradient of the loss function with respect to $\x{B}$, thus requiring $O(KN^3)$ time per iteration.

Finally, that approach does not utilize further computational speed-ups that are possible by enforcing a simple orthogonality constraint on $\x{B}$, \textit{viz}., where we require $\x{B}$ to be orthonormal (i.e., $\x{B}\xt{B} = \x{I}_N$). Here we note that although this so-called orthogonal joint diagonalization (OJD) property \cite{mesloub2014} is built into several existing methods, most notably JADE \cite{cardoso1996}, those methods still use sweeps, where each sweep also requires $O(KN^3)$ time.

Instead, we (\emph{i}) follow a quasi-Newton approach to jointly update all elements $\x{B}$ in each iteration, (\emph{ii}) exploit orthonormality to further simplify the criterion, (\emph{iii}) parametrize the optimization problem such that we can use unconstrained optimization while guaranteeing orthonormality of $\x{B}$, (\emph{iv}) use a low-dimensional approximation of $\x{C}_k$ to avoid calculation of matrix products of $N \times N$ matrices for $k=1,\ldots,K$ in each iteration, (\emph{v}) apply a mild regularization to the approximations of matrices $\x{C}_k$ such that it does not perturb our solution for $\x{B}$, and (\emph{vi}) linearize updates to efficiently implement a golden-section search \cite{kiefer1953} within each iteration.

This overall approach, incorporating these six points, we refer to as JADOC: joint approximate diagonalization under orthogonality constraints. JADOC requires $O(N^3)$ time per iteration, constituting an $O(K)$ reduction in runtime compared to existing methods. Moreover, the marginal cost of performing the golden-section search within each iteration is negligible.

Our main derivations are provided in Section \ref{sec:main}, our JADOC algorithm is presented in Section \ref{sec:alg} together with the algorithm that we use to simulate data, simulation results are discussed in Section \ref{sec:results}, and the conclusions follow in Section \ref{sec:conclusions}.

\section{Main derivations}\label{sec:main}
We follow Pham's notation \cite{pham2001}, while using the criterion as proposed by Ablin et al. \cite{ablin2018}, which is defined as
\begin{equation}\label{eq:crit_ablin}
\mathcal{L}\left( \x{B} \right) = \frac{1}{2K} \sum_{k=1}^K \left(\log \det \diag \left( \x{B} \x{C}_k \xt{B} \right) - \log \det \left( \x{B} \x{C}_k \xt{B} \right) \right).
\end{equation}
This criterion is in fact equivalent to that in Equation \ref{eq:crit_pham}, except for the scaling by $n_k$, which effectively weights $\x{C}_k$ by its degrees of freedom. We ignore this scaling, as we do not necessarily assume $\x{C}_k$ to follow a Wishart distribution. We simply want to find $\x{B}$ that diagonalizes symmetric, PSD matrices $\x{C}_k$ for $k=1,\ldots,K$ as much as possible, using the criterion in Equation \ref{eq:crit_ablin} (i.e., without assigning particular weights to the individual matrices $\x{C}_k$).

Given that (\emph{i}) $\det (\x{A}_1\x{A}_2 ) = \det (\x{A}_1)\det (\x{A}_2)$ for square matrices $\x{A}_1$ and $\x{A}_2$, and (\emph{ii}) $\x{B}\xt{B} = \x{I}_N$, we can rewrite this criterion as follows
\begin{align*}
\mathcal{L}\left( \x{B} \right)  &= \frac{1}{2K} \sum_{k=1}^K \left(\log \det \diag \left( \x{B} \x{C}_k \xt{B} \right) - \log \left( \det \left( \x{B} \right) \det \left( \x{C}_k \right) \det \left( \xt{B} \right) \right) \right) \\
&= \frac{1}{2K} \sum_{k=1}^K \left(\log \det \diag \left( \x{B} \x{C}_k \xt{B} \right) - \log \left(  \det \left( \x{C}_k \right) \det \left( \x{B} \xt{B} \right) \right) \right) \\
&= \frac{1}{2K} \sum_{k=1}^K \left(\log \det \diag \left( \x{B} \x{C}_k \xt{B} \right) - \log \left(  \det \left( \x{C}_k \right)  \right) \right).
\end{align*}
Clearly, the second log-determinant in the summand is independent of $\x{B}$. Consequently, this term merely contributes to the constant. Thus, we can safely drop the second term from the loss function. Moreover, let the eigenvalue decomposition of $\x{C}_k$ be given by
\begin{equation}\label{eq:evd}
\x{C}_k = \x{P}_k \x{D}_k^2 \xt{P}_k,
\end{equation}
where $\x{P}_k$ denotes the orthonormal matrix of eigenvectors and $\x{D}_k^2$ the diagonal matrix of eigenvalues. Now, further, consider a rank-$S$ approximation to $\x{C}_k$, defined as
\begin{equation}\label{eq:rank_S_approx}
\widetilde{\x{C}}_k = \x{L}_k \xt{L}_k, \textrm{ where } \x{L}_k = \widetilde{\x{P}}_k \widetilde{\x{D}}_k ,
\end{equation}
where $\widetilde{\x{D}}_k$ is an $S \times S$ diagonal matrix, with the square root of the $S$ leading eigenvalues of $\x{C}_k$ as its diagonal elements, and $\widetilde{\x{P}}_k$ denotes the $N \times S$ matrix of corresponding eigenvectors.

JADOC aims to find a matrix $\x{B}$ that diagonalizes this rank-$S$ approximation, $\widetilde{\x{C}}_k$, for $k=1,\ldots,K$. Importantly, when $\widetilde{\x{C}}_k$ is rank-deficient, for at least one $k \in \{1,\ldots,K\}$, the loss function is ill-defined. However, observe that a matrix $\x{B}$ that diagonalizes $\widetilde{\x{C}}_k$ also diagonalizes $\widetilde{\x{C}}_k + \lambda \x{I}$ for some fixed value of $\lambda$.

Moreover, by setting $\lambda > 0$, matrix $\widetilde{\x{C}}_k + \lambda \x{I}$ is positive definite for $k=1,\ldots,K$. Diagonalization of these matrices is, thus, a valid approach to find a matrix $\x{B}$ that also (approximately) diagonalizes $\widetilde{\x{C}}_k$, with the additional advantage that the loss function is properly defined for all orthonormal real matrices $\x{B}$. Thus, more concretely, the aim of JADOC is to find an orthonormal, real matrix $\x{B}$ that diagonalizes the rank-$S$ approximations of $\x{C}_k$ for $k=1,\ldots,K$, by minimizing:
\begin{align}
\mathcal{L}\left( \x{B} \right)  &= \frac{1}{2K} \sum_{k=1}^K  \log \det \diag \left( \x{B} \left( \widetilde{\x{C}}_k + \lambda \x{I} \right)  \xt{B} \right), \textrm{ where} \\
\lambda &= \lambda_0 + \frac{1}{NK} \sum_{k=1}^K \tr \left( \x{C}_k - \widetilde{\x{C}}_k \right), \textrm{ where } \lambda_0 = 1. \label{eq:lambda}
\end{align}
The first term in our definition of $\lambda$ (i.e., $\lambda_0$) ensures the loss function is well-defined even though $\widetilde{\x{C}}_k$ is rank deficient for sure when $S < N$ and irrespective of whether the original input matrices $\x{C}_k$ are rank deficient or not. The second term ensures the scale of the diagonal elements in $\widetilde{\x{C}}_k + \lambda \x{I}$ is on average roughly the same as the diagonal elements of $\x{C}_k$. Further, we note that the trace in Equation \ref{eq:lambda} simply equals the trace of $\x{C}_k$ minus the sum of its $S$ leading eigenvalues.

We can now simplify the preceding criterion to the \textbf{main JADOC criterion}:
\begin{equation}\label{eq:crit_devlaming}
\mathcal{L} \left( \x{B} \right) = \frac{1}{2K} \sum_{k=1}^K \sum_{i=1}^N \log \left( \lambda + \sum_{j=1}^S \left( \x{A}_k \right)_{ij}^2 \right), \textrm{ where } \x{A}_k = \x{B} \x{L}_k,
\end{equation}
where $\lambda$ is as defined in Equation \ref{eq:lambda}. Now, by setting
\begin{equation}
S = \left\lceil \frac{N}{K} \right\rceil,\label{eq:s}
\end{equation}
and given (\emph{i}) we pre-calculate $\x{L}_k$, requiring $O(KN^3)$ time in total, and (\emph{ii}) we want to evaluate the loss function for some orthonormal $\x{B}$, we can calculate $\x{A}_k$ for $k=1,\ldots,K$ jointly in $O(N^3)$ time. This time complexity is the bottleneck for calculations in each iteration.

We consider the choice of $S$ in Equation \ref{eq:s} to be fair: (\emph{i}) it increases with $N$, (\emph{ii}) it extracts the most salient dimensions from each matrix $\x{C}_k$, for $k=1,\ldots,K$, and (\emph{iii}) it is the lowest $S$ such that there may exist convex combinations of $\widetilde{\x{C}}_k$, for $k=1,\ldots,K$, that are positive definite. Notice that users of the JADOC tool may specify $S$ different from the recommended value seen in Equation \ref{eq:s}.

\subsection{Parametrization}
Given our loss function in Equation \ref{eq:crit_devlaming} and the fact that we will use an iterative procedure to obtain $\x{B}$, let us assume that in iteration $t$ we have $\x{B}$ from the directly preceding iteration, denoted by $\x{B}_{(t-1)}$. Here, we further assume that $\x{B}_{(t-1)}$ is orthonormal. That is, $\x{B}_{(t-1)} \xt{B}_{(t-1)} = \x{I}_N$.

Bearing these considerations in mind, the goal of iteration $t$ is to `update' $\x{B}_{(t-1)}$ such that (\emph{i}) the loss function is further decreased and (\emph{ii}) orthonormality of $\x{B}$ is preserved. The easiest way to ensure orthogonality is indeed maintained, is by using a rotation matrix $\x{R}$ (i.e., a matrix such that $\x{R}\xt{R} = \x{I}_N$) to set $\x{B}_{(t)} = \x{R} \x{B}_{(t-1)}$, as clearly under this definition of $\x{B}_{(t)}$, we have that $\x{B}_{(t)}\xt{B}_{(t)} = \x{I}_N$.

However, defining such an $N \times N$ matrix $\x{R}$ `straight up' is far from trivial. Fortunately, for any rotation matrix $\x{R}$ there exists a skew-symmetric matrix $\x{S}$ such that $e^{\x{S}} = \x{R}$ \cite{gallier2002,brocker1985}, where $e^{c\x{M}}$ denotes the matrix exponential of square matrix $\x{M}$, which is defined in terms of the following power series \cite{marsden1999}:
\begin{equation}\label{eq:power}
e^{c\x{M}} = \x{I} + \sum_{k=1}^{\infty} \frac{c^k }{k!} \x{M}^k.
\end{equation}
Thus, we can achieve our second aim in any iteration by defining a strictly lower-triangular matrix $\x{E}$ (i.e., with diagonal elements equal to zero), and considering each of its free elements $\varepsilon_{lm} = \left( \x{E} \right)_{lm} \in \mathbb{R}$ for $l,m=1,\ldots,N$ where $l>m$ as a parameter of interest. Under this parametrization, (\emph{i}) $\x{S} = \x{E} - \xt{E}$ defines a skew-symmetric matrix (i.e., such that $\xt{S} = -\x{S}$), for which it holds that $e^{\x{S}}$ is a rotation matrix $\x{R}$ (i.e., such that $\x{R}\xt{R} = \x{I}_N$), and (\emph{ii}) any $\x{B}$ that meets the orthonormality requirement can be written as follows
\begin{equation}
\x{B}_{(t)} = e^{\x{E} - \xt{E}} \x{B}_{(t-1)}.
\end{equation}
Under these definitions, letting $\x{A}_{(t)k}$ denote $\x{A}_k$ in iteration $(t)$, notice that
\begin{equation}
\x{A}_{(t)k} = \x{B}_{(t)} \x{L}_k = e^{\x{E} - \xt{E}}\x{B}_{(t-1)}\x{L}_k = e^{\x{E} - \xt{E}} \x{A}_{(t-1)k},
\end{equation}
implying we can define both $\x{B}$ and $\x{A}_k$ in a recursive manner across iterations.

Bearing all these considerations in mind, we can write our loss function in terms of the free parameters in $\x{E}$ and $\x{A}_k$ from the preceding iteration as follows:
\begin{align}\label{eq:crit_final}
\mathcal{L} \left( \x{E} \right) &= \frac{1}{2K} \sum_{k=1}^K \sum_{i=1}^N \log \left( d_{ik} \left( \x{E} \right) \right), \textrm{ where } \\
d_{ik} \left( \x{E} \right) &= \lambda + \sum_{j=1}^S \left( e^{\x{E} - \xt{E}} \x{A}_{(t-1)k} \right)_{ij}^2.
\end{align}

\subsection{First derivatives}
Here, let us consider the partial derivative of this loss function with respect to $\varepsilon_{lm}$ for $l>m$. This derivative can be written as
\begin{equation}\label{eq:grad1}
\frac{\partial \mathcal{L} \left( \x{E} \right)}{\partial \varepsilon_{lm}}  = \frac{1}{K} \sum_{k=1}^K \sum_{i=1}^N \frac{1}{ d_{ik} \left( \x{E} \right) } \sum_{j=1}^N \left( e^{\x{E} - \xt{E}} \x{A}_{(t-1)k} \right)_{ij} \frac{\partial \left( e^{\x{E} - \xt{E}} \x{A}_{(t-1)k} \right)_{ij} }{\partial \varepsilon_{lm}},
\end{equation}
where
\begin{equation*}
\left( e^{\x{E} - \xt{E}} \x{A}_{(t-1)k} \right)_{ij} = \sum_{p=1}^N \left( e^{\x{E} - \xt{E}} \right)_{ip} \left( \x{A}_{(t-1)k} \right)_{pj}.
\end{equation*}
Furthermore, using the power series in Equation \ref{eq:power}, it is easy to show that
\begin{equation*}
\frac{\partial e^{\x{E} - \xt{E}} }{\partial \varepsilon_{lm}} = e^{\x{E} - \xt{E}} \left( \x{J}_{lm} - \x{J}_{ml} \right),
\end{equation*}
where $\x{J}_{lm}$ is the single-entry matrix, which equals one for element $l,m$ and zero elsewhere. Thus, we have that
\begin{align*}
\frac{\partial \left( e^{\x{E} - \xt{E}} \x{A}_{(t-1)k} \right)_{ij} }{\partial \varepsilon_{lm}} &= \sum_{p=1}^N \frac{\partial \left( e^{\x{E} - \xt{E}} \right)_{ip} }{\partial \varepsilon_{lm}} \left( \x{A}_{(t-1)k} \right)_{pj} \\
&= \sum_{p=1}^N \left( e^{\x{E} - \xt{E}} \left( \x{J}_{lm} - \x{J}_{ml} \right) \right)_{ip}  \left( \x{A}_{(t-1)k} \right)_{pj}
\end{align*}
Substituting this expression in Equation \ref{eq:grad1} yields the following expression for the partial derivative of the loss function with respect to $\varepsilon_{lm}$:
\begin{align*}
\frac{\partial \mathcal{L} \left( \x{E} \right)}{\partial \varepsilon_{lm}} &= \frac{1}{K} \sum_{k=1}^K \sum_{i=1}^N \frac{1}{ d_{ik}\left( \x{E} \right)} \sum_{j=1}^N \left[ \left( e^{\x{E} - \xt{E}} \x{A}_{(t-1)k} \right)_{ij} \right. \\
& \hphantom{ = \frac{1}{K} } \left. \sum_{p=1}^N \left( e^{\x{E} - \xt{E}} \left( \x{J}_{lm} - \x{J}_{ml} \right) \right)_{ip}  \left( \x{A}_{(t-1)k}  \right)_{pj} \right] \\
&= \frac{1}{K} \sum_{k=1}^K \sum_{i=1}^N \frac{1}{ d_{ik} \left( \x{E} \right)} \left( e^{\x{E} - \xt{E}} \x{A}_{(t-1)k} \xt{A}_{(t-1)k} \left( \x{J}_{ml} - \x{J}_{lm}\right) e^{\xt{E} - \x{E}} \right)_{ii}
\end{align*}
Evaluation of the loss function and gradient at current value of $\x{B}$, which corresponds to $\x{E} = \x{0}$ and $e^{\x{E}-\xt{E}} = \x{I}_N$, yields
\begin{align}\label{eq:grad_eval}
{\left. \frac{\partial \mathcal{L} \left( \x{E} \right)}{\partial \varepsilon_{lm}} \right\vert}_{\x{E}=\x{0}} &= \frac{1}{K} \sum_{k=1}^K \sum_{i=1}^N \frac{1}{ d_{ik}} \left( \x{A}_{(t-1)k} \xt{A}_{(t-1)k} \left( \x{J}_{ml} - \x{J}_{lm}\right) \right)_{ii}, \textrm{ where } \\
d_{ik} &= d_{ik} \left( \x{0} \right) \equiv \lambda + \sum_{j=1}^S \left( \x{A}_{(t-1)k} \right)^2_{ij}
\end{align}
Using properties of single-entry matrices, this expression can be simplified as follows:
\begin{equation}\label{eq:grad_eval_simple}
{\left. \frac{\partial \mathcal{L} \left( \x{E} \right)}{\partial \varepsilon_{lm}} \right\vert}_{\x{E}=\x{0}} = \frac{1}{K} \sum_{k=1}^K \left[ \left( \frac{1}{ d_{lk}} - \frac{1}{ d_{mk}} \right) \left( \x{A}_{(t-1)k} \xt{A}_{(t-1)k}  \right)_{lm}  \right].
\end{equation}
Now, the gradient can be expressed in matrix notation as follows:
\begin{align}
\x{G} &= \frac{\partial \mathcal{L} \left( \x{E} \right) }{\partial \x{E}} = \stril \left( \x{F} - \xt{F} \right), \textrm{ where} \label{eq:grad_matrix} \\
\x{F} &= \frac{1}{K} \sum_{k=1}^K \diag \left(  \frac{1}{d_{1k}}, \ldots,  \frac{1}{d_{Nk}} \right) \left( \x{A}_{(t-1)k} \xt{A}_{(t-1)k} \right),
\end{align}
where $\stril (\x{F} - \xt{F})$ is the strict lower-triangular submatrix of square matrix $\x{F} - \xt{F}$.

\subsection{Second derivatives}
Let us now consider second partial derivatives, which are, in general, defined by the following expression for $l>m$ and $p>q$:
\begin{align}
\frac{\partial^2 \mathcal{L} \left( \x{E} \right) }{\partial \varepsilon_{lm} \partial \varepsilon_{pq}} &= \frac{1}{K} \sum_{k=1}^K \sum_{i=1}^N \frac{ \partial \left( \frac{1}{ d_{ik} \left( \x{E} \right)} \left( \x{T}_k \left( \x{E} \right) \right)_{ii} \right)}{\partial \varepsilon_{pq}}, \textrm{ where} \label{eq:2nd} \\
\x{T}_k \left( \x{E} \right) &= e^{\x{E} - \xt{E}} \x{A}_{(t-1)k} \xt{A}_{(t-1)k} \left( \x{J}_{ml} - \x{J}_{lm}\right) e^{\xt{E} - \x{E}}.
\end{align}
Using the chain rule, the partial derivative on the right-hand side of Equation \ref{eq:2nd} can be written as follows:
\begin{equation}\label{eq:2nd_inner}
\frac{ \partial \left( \frac{1}{ d_{ik} \left( \x{E} \right)} \left( \x{T}_k \left( \x{E} \right) \right)_{ii} \right)}{\partial \varepsilon_{pq}} = \frac{1}{ d_{ik}\left( \x{E} \right)} \frac{\partial \left( \x{T}_k \left( \x{E} \right) \right)_{ii}}{\partial \varepsilon_{pq}} - \frac{1}{\left( d_{ik}\left( \x{E} \right)\right)^2} \frac{\partial d_{ik}\left( \x{E} \right)}{\partial \varepsilon_{pq}}  \left( \x{T}_k \left( \x{E} \right) \right)_{ii}.
\end{equation}
Immediately evaluating $\left( \x{T}_k \left( \x{E} \right) \right)_{ii}$ in the second term at $\x{E}=\x{0}$, yields
\begin{equation}
 \left( \x{T}_k \left( \x{0} \right) \right)_{ii} = \left\{ \begin{array}{ll}\left( \x{A}_{(t-1)k} \xt{A}_{(t-1)k}  \right)_{lm} & \textrm{ if } i=l,  \\ - \left( \x{A}_{(t-1)k}  \xt{A}_{(t-1)k}  \right)_{lm} & \textrm{ if } i=m, \\ 0 & \textrm{ otherwise. } \end{array}\right.
\end{equation}
Effectively, this implies the second term on the right-hand side of Equation \ref{eq:2nd_inner} is zero when $\x{A}_{(t-1)k}  \xt{A}_{(t-1)k} $ is diagonal. Thus, assuming (near) diagonalization is possible and has already been (almost) achieved in iteration $t-1$, this last term can be ignored.

Consequently, the most important part of the second partial derivatives is finding an efficient expression for
\begin{equation}
\frac{\partial \left( \x{T}_k \left( \x{E} \right) \right)_{ii}}{\partial \varepsilon_{pq}} = \frac{\partial \left( e^{\x{E} - \xt{E}} \x{A}_{(t-1)k}  \xt{A}_{(t-1)k}  \left( \x{J}_{ml} - \x{J}_{lm}\right) e^{\xt{E} - \x{E}} \right)_{ii} }{\partial \varepsilon_{pq}}.
\end{equation}
We here, thus, assume that $\x{A}_{(t-1)k}  \xt{A}_{(t-1)k} $ is diagonal. That is, we can replace $\x{A}_{(t-1)k}  \xt{A}_{(t-1)k} $ by diagonal matrix $\x{D}_k$ of which element $i,i$ is given by $d_{ik} - \lambda$. Therefore, we have that
\begin{equation}
\frac{\partial \left( \x{T}_k \left( \x{E} \right) \right)_{ii}}{\partial \varepsilon_{pq}} = \sum_{s=1}^N \sum_{t=1}^N \frac{ \partial \left(  \left( e^{\x{E} - \xt{E}} \right)_{it} \left( \x{D}_k \left( \x{J}_{ml} - \x{J}_{lm}\right)\right)_{ts} \left( e^{\xt{E} - \x{E}} \right)_{si} \right) }{\partial \varepsilon_{pq}}.
\end{equation}
Using the product rule and previously discussed properties of the partial derivative of the matrix exponential with respect to a skew symmetric matrix, it then follows that
\begin{align*}
\frac{\partial \left( \x{T}_k \left( \x{E} \right) \right)_{ii}}{\partial \varepsilon_{pq}} &= \sum_{s=1}^N \sum_{t=1}^N \left[ \frac{ \partial \left( e^{\x{E} - \xt{E}} \right)_{it} }{\partial \varepsilon_{pq}} \left( \x{D}_k \left( \x{J}_{ml} - \x{J}_{lm}\right)\right)_{ts} \left( e^{\xt{E} - \x{E}} \right)_{si} \right. \\
& \hphantom{=} \left. + \left( e^{\x{E} - \xt{E}} \right)_{it} \left( \x{D}_k \left( \x{J}_{ml} - \x{J}_{lm}\right)\right)_{ts}  \frac{ \partial \left(  e^{\xt{E} - \x{E}} \right)_{si}}{\partial \varepsilon_{pq}} \right] \\
&= \sum_{s=1}^N \sum_{t=1}^N \left[ \left( \left( e^{\x{E} - \xt{E}} \right) \left( \x{J}_{pq} - \x{J}_{qp}\right) \right)_{it}  \left( \x{D}_k \left( \x{J}_{ml} - \x{J}_{lm}\right)\right)_{ts} \left( e^{\xt{E} - \x{E}} \right)_{si} \right. \\
& \hphantom{=} \left. + \left( e^{\x{E} - \xt{E}} \right)_{it} \left( \x{D}_k \left( \x{J}_{ml} - \x{J}_{lm}\right)\right)_{ts} \left(  \left( \x{J}_{qp} - \x{J}_{pq}\right) \left( e^{\xt{E} - \x{E}} \right) \right)_{si} \right] \\
&= \left( \left( e^{\x{E} - \xt{E}} \right) \left( \x{J}_{pq} - \x{J}_{qp}\right)  \x{D}_k \left( \x{J}_{ml} - \x{J}_{lm}\right) \left( e^{\xt{E} - \x{E}} \right) \right)_{ii} \\
& \hphantom{=} + \left( \left( e^{\x{E} - \xt{E}} \right) \x{D}_k \left( \x{J}_{ml} - \x{J}_{lm}\right) \left( \x{J}_{qp} - \x{J}_{pq}\right) \left( e^{\xt{E} - \x{E}} \right) \right)_{ii}.
\end{align*}
Evaluating this expression at $\x{E} = \x{0}$ and substituting in our expression for the second-order derivative under complete diagonality, we have that
\begin{align*}
\left. \frac{\partial^2 \mathcal{L} \left( \x{E} \right) }{\partial \varepsilon_{lm} \partial \varepsilon_{pq}} \right\vert_{\x{E} = \x{0}} &= \frac{1}{K} \sum_{k=1}^K \sum_{i=1}^N \frac{1}{ d_{ik}} \left[ \left( \left( \x{J}_{pq} - \x{J}_{qp}\right)  \x{D}_k \left( \x{J}_{ml} - \x{J}_{lm}\right)  \right)_{ii} \right. \\
& \hphantom{= \frac{1}{K}} \left. + \left(  \x{D}_k \left( \x{J}_{ml} - \x{J}_{lm}\right) \left( \x{J}_{qp} - \x{J}_{pq}\right) \right)_{ii} \right].
\end{align*}
The second term in the inner summand is easiest to simplify further. Notice that by properties of the single-entry matrix, this term is zero, except when $l=p$ and $m=q$. That is, when we consider the second-order derivative of the loss function with respect to $\varepsilon_{lm}$. In that particular case, the second term is given by $- \left( \x{D}_k \x{J}_{mm} + \x{D}_k \x{J}_{ll} \right)_{ii}$. For the first term, things are slightly more involved. Yet, here, we can show that it equals zero, again except when $l=p$ and $m=q$, in which case it equals the following expression: $\left( \left( d_{mk} - \lambda \right) \x{J}_{ll} + \left( d_{lk} - \lambda \right) \x{J}_{mm} \right)_{ii}$.

Now, considering index $i$, for all values of $i$ other than $m$ and $l$, both terms are also zero, even if $l=p$ and $m=q$. Consequently, we can write the second-order derivative with respect to $\varepsilon_{lm}$ as follows:
\begin{align}
\left. \frac{\partial^2 \mathcal{L} \left( \x{E} \right) }{\partial^2 \varepsilon_{lm}} \right\vert_{\x{E} = \x{0}} &= \frac{1}{K} \sum_{k=1}^K \sum_{i=1}^N \frac{1}{ d_{ik}} \left[ \left( \left( d_{mk} - \lambda \right) \x{J}_{ll} + \left( d_{lk} - \lambda \right) \x{J}_{mm}  \right)_{ii} \right. \\
& \hphantom{= \frac{1}{K}} \left. - \left(  \x{D}_k \x{J}_{mm} + \x{D}_k \x{J}_{ll} \right)_{ii} \right] \\
&= \frac{1}{K} \sum_{k=1}^K \frac{ \left( d_{mk} - \lambda \right)  - \left( d_{lk}  - \lambda \right) }{ d_{lk}} + \frac{ \left( d_{lk} - \lambda \right) - \left( d_{mk} - \lambda \right) }{ d_{mk}} \\
&= \frac{1}{K} \sum_{k=1}^K \frac{d_{mk}}{d_{lk}} + \frac{d_{lk}}{d_{mk}} - 2. \label{eq:hessian_diag}
\end{align}
In short, when the solution is very close to being diagonal, the Hessian becomes a diagonal matrix, with the diagonal element that corresponds to the second derivative with respect to parameter $\varepsilon_{lm}$ as described in the last equation. Thus, by dividing the corresponding elements in the gradient by these elements, we obtain a quasi-Newton update.

Defining matrix $\x{H}$ with element $l,m$ as seen in Equation \ref{eq:hessian_diag}, the search direction in this minimization problem can now be set as follows:
\begin{equation}
\x{E} = - \x{G} \oslash \x{H}, \label{eq:update}
\end{equation}
where `$\oslash$' denotes the elementwise division, of $N \times N$ matrices $\x{G}$ and $\x{H}$, comprising corresponding elements of the gradient and Hessian respectively.

Each element of $\x{H}$ can be described as a function of some scalar $x \in [0,+\infty)$, \textit{viz}., $x + x^{-1} -2$. It is easy to show this function attains its minimum at $x=1$, where then the value is exactly zero. Thus, the elements of $\x{H}$ are non-negative by definition. To ensure numerical stability, we set elements of $\x{H}$ below threshold $\tau_H = 0.01$ to be equal to $\tau_H$. This modification to $\x{H}$ is equivalent to requiring that some element in $\x{E}$ can have a magnitude that is at most 100 times as large as the corresponding element in $\x{G}$.

This overall approach for calculating $\x{G}$, $\x{H}$, and $\x{E}$ is similar to what is used in an undocumented part of the Python package \texttt{qndiag} \cite{ablin2018}. However, our derivations (\emph{i}) formalize and prove the correctness of expressions that are similar to those utilized in that code, (\emph{ii}) provide expressions that are formulated such that greater computational speedups are possible, and (\emph{iii}) take our low-rank approximation with regularization fully into account.

\subsection{Line search}
\label{subsec:line_search}
The preceding quasi-Newton approach formulates the update in terms of matrix $\x{E}$. In a standard line search, one considers $\alpha \in [0,1]$, where $\alpha=0$ would correspond to no update and $\alpha=1$ would correspond to the full update. In our formulation, the update has a multiplicative nature.

That is, define
\begin{equation}
\x{R}^* = e^{\x{E} - \xt{E}}.
\end{equation}
Now, considering $\x{A}_{(t-1)k}$ for $k=1,\ldots,K$ (readily available from the previous iteration), if $\alpha=1$, then $\x{A}_{(t)k}$ would be given by $\x{R}^* \x{A}_{(t-1)k}$. Thus, given $\x{E}$ and $\x{A}_{(t-1)k}$, we have that $\x{R}^*$ and $\x{R}^* \x{A}_{(t-1)1},\ldots,\x{R}^* \x{A}_{(t-1)K}$ altogether can be calculated in $O(N^3)$ time. Moreover, if $\alpha=0$, we have that $\x{A}_{(t)k}$ is simply given by its value in the previous iteration.

By taking a convex combination of these two end-points, we `linearize' the update during the line search. That is, we set
\begin{equation}
\x{A}^*_k = \alpha \x{R}^* \x{A}_{(t-1)k} + (1- \alpha) \x{A}_{(t-1)k}.
\end{equation}
Now, given $\x{R}^* \x{A}_{(t-1)k}$ and $\x{A}_{(t-1)k}$ for $k=1,\ldots,K$, we can calculate the loss function as defined in Equation \ref{eq:crit_devlaming} in $O(KN^2)$ time for each value of $\alpha$ considered. This approach allows us to efficiently apply a golden-section search to find the best value of $\alpha$. 

Finally, we observe that as a result of the linearization, for $\alpha \in (0,1)$ its value is slightly underestimated. Therefore, given $\alpha$ after the linearization, we cast it back to the exponential scale by setting
\begin{equation}
\alpha^* = \log \left(1+\left(\alpha(e-1)\right) \right).
\end{equation}
Finally, we set 
\begin{align}
\x{R} &= e^{\alpha^* \left( \x{E} - \xt{E}\right) } \label{eq:final_rotation}\\
\x{B}_{(t)} &= \x{R} \x{B}_{(t-1)} \label{eq:update_b} \\
\x{A}_{(t)k} &= \x{R} \x{A}_{(t-1)k} \label{eq:update_a},
\end{align}
after which we move to the next iteration, provided convergence has not yet occurred.

\subsection{Initialization}
JADOC initializes by setting $\x{B}_{(0)} = \x{I}_N$, which clearly is a valid orthonormal transformation matrix.

\subsection{Convergence}
Convergence is defined in terms of the root-means-square deviation (RMSD) from zero of the gradient. By default, a value below $10^{-4}$, in conjunction with JADOC having made at least $T_{min}= 10$ iterations, is considered to constitute convergence. The latter requirement is to deal with the fact that the gradient can be quite small near the starting point (i.e., in the first few iterations) due to poor starting values. After $T=100$ iterations JADOC will terminate, irrespective of whether convergence has occurred.

\section{Algorithms}
\label{sec:alg}
Algorithm \ref{alg:jadoc} shows an overview of the JADOC algorithm as implemented in our Python code, available at \url{https://github.com/devlaming/jadoc}. Within the iterations (i.e., the last for loop), updating $\x{B}$ and $\x{A}_k$ for $k=1,\ldots,K$ (i.e., the last step) is the most demanding step, computationally speaking, causing the algorithm to require $O(N^3)$ per iteration. This, however, still constitutes an $O(K)$ reduction in time complexity, compared to existing algorithms for joint diagonalization.

Algorithm \ref{alg:simulation} shows an overview of the data-generating process we use in our simulations, to compare performance of JADOC to existing algorithms. This algorithm is implemented in the \texttt{SimulateData(iK,iN,iR,dAlpha)} function, as implemented in the Python package on GitHub.

\begin{algorithm}
\caption{JADOC algorithm.}
\label{alg:jadoc}
\hspace*{\algorithmicindent}\textbf{Input}: $N\times N$ positive (semi)-definite matrices, $\x{C}_k$, for $k=1,\ldots,K$ \\
\hspace*{\algorithmicindent}\textbf{Output}: $N\times N$ orthonormal diagonalizer matrix
\begin{algorithmic}
\STATE{Set $\tau_{conv} = 10^{-4}$, $\tau_{H} = 10^{-2}$, $S = \left\lceil \frac{N}{K} \right\rceil$, $\lambda_0 = 1$, $\x{B}_{(0)} = \x{I}_N$, $T=100$, $T_{min}=10$}
\FOR{$k=1,\ldots,K$}
\STATE{Compute $\x{L}_k$ in line with Equation \ref{eq:rank_S_approx} and store as $\x{A}_{(0)k}$}
\STATE{Update $\lambda$ in line with Equation \ref{eq:lambda}}
\ENDFOR
\FOR{$t=1,\ldots,T$}
\STATE{Compute gradient in accordance with Equation \ref{eq:grad_matrix} $\rightarrow \x{G}$}
\IF{root-mean-square deviation of gradient $< \tau_{conv}$ \textbf{and} $t > T_{min}$}
\STATE{\textbf{break for loop}}
\ENDIF
\STATE{Compute diagonal of approximated Hessian using Equation \ref{eq:hessian_diag} $\rightarrow \x{H}$}
\STATE{Set elements of $\x{H}$ below $\tau_H$ equal to $\tau_H$}
\STATE{Compute update using Equation \ref{eq:update} $\rightarrow \x{E}$}
\STATE{Apply line search described in Subsection \ref{subsec:line_search} $\rightarrow \alpha^*$}
\STATE{Computation rotation matrix using Equation \ref{eq:final_rotation} $\rightarrow \x{R}$}
\STATE{Update $\x{B}$ and $\x{A}_k$ using Equations \ref{eq:update_b} \& \ref{eq:update_a} $\rightarrow \x{B}_{(t)}, \x{A}_{(t)k}$}
\ENDFOR
\RETURN $\x{B}_{(t)}$
\end{algorithmic}
\end{algorithm}

\begin{algorithm}
\caption{Algorithm for the data-generating process of our simulations.}
\label{alg:simulation}
\hspace*{\algorithmicindent}\textbf{Input}: $K$, $N$, replicate number $r$, degree of similarity eigenvectors $\alpha \in [0,1]$\\
\hspace*{\algorithmicindent}\textbf{Output}: $K$ different $N\times N$ PSD matrices
\begin{algorithmic}
\STATE{Set random seed based on $r$}
\STATE{Get $N^2$ iid draws from $\mathcal{N} \left( 0,1\right)$ and store in $N \times N$ matrix $\rightarrow \x{X}$}
\FOR{$k=1,\ldots,K$}
\STATE{Get $N^2$ iid draws from $\mathcal{N} \left( 0,1\right)$ and store in $N \times N$ matrix $\rightarrow \x{X}_k $}
\STATE{Compute $\alpha \x{X} + \left(1 - \alpha \right) \x{X}_k \rightarrow \x{X}_k$}
\STATE{Compute $e^{\x{X}_k - \xt{X}_k} \rightarrow \x{R}$}
\STATE{Get matrix of zeros $\rightarrow \x{D}$}
\FOR{$n=1,\ldots,N$}
\STATE{Draw from $\chi^2 \left(1 \right) \rightarrow d$}
\STATE{Set $n$-th diagonal element of $\x{D}$ equal to $d$}
\ENDFOR
\STATE{Compute $\x{R} \x{D} \xt{R} \rightarrow \x{C}_k$}
\ENDFOR
\RETURN $\x{C}_{1},\ldots,\x{C}_K$
\end{algorithmic}
\end{algorithm}

\section{Simulation results}
\label{sec:results}

We consider two simulation designs. In Design 1, we set set $K=10$ and consider $N=100, 200, 300, 400, 500$. In Design 2, we set $N=256$ and consider $K=2,4,8,16,32$. In each design, for a given value of $N$ and $K$, we consider $R=10$ replicates.

For each of these designs, we consider four values of the parameter $\alpha$ defined in Algorithm \ref{alg:simulation}, \textit{viz}., $\alpha=0,0.25,0.5,0.75$. Here, $\alpha=0$ corresponds to $\x{C}_1, \ldots, \x{C}_K$ having completely independent eigenvectors. Values of $\alpha$ between zero and one correspond to intermediate scenarios, where eigenvectors across matrices are correlated; the higher $\alpha$, the stronger this correlation.

In total, we have eight scenarios: two designs and four different values of $\alpha$. In each scenario, we apply \begin{enumerate}
\item JADE \cite{cardoso1996} as implemented on \\ {\scriptsize{\url{https://github.com/gabrieldernbach/approximate_joint_diagonalization}}},
\item \texttt{qndiag} \cite{ablin2018} as implemented on \\ {\scriptsize{\url{https://github.com/pierreablin/qndiag}}}, and
\item our JADOC algorithm as implemented on \\ {\scriptsize{\url{https://github.com/devlaming/jadoc}}}.
\end{enumerate}

For each of the ten replicates and each of the three methods, we run the relevant Python code on separate machines with identical configurations (i.e., each machine has 24 cores with a clock speed of 2.6GHz and 64GB of memory). That is, for a given replicate and a given method, all ten scenarios are carried out on a single, separate machine. For JADOC and JADE, we make use of the built-in parallelized functionality. Finally, for fair comparison, for \texttt{qndiag}, we use the built-in option to enforce orthonormality of $\x{B}$.

For each scenario and for each method, we report the median (across the ten replicates) of the runtimes (solid lines) and the median (across the ten replicates) of the root-mean-square deviation from zero of the off-diagonal elements after transformation using the matrix $\x{B}$ (dotted lines).
\begin{figure}[p!]
  \centering
  \begin{tabular}{cc} \hline
  \multicolumn{2}{c}{$\alpha=0.00$} \\
  \includegraphics[scale=0.4,trim=0.7in 0in 0.7in 0.35in,clip]{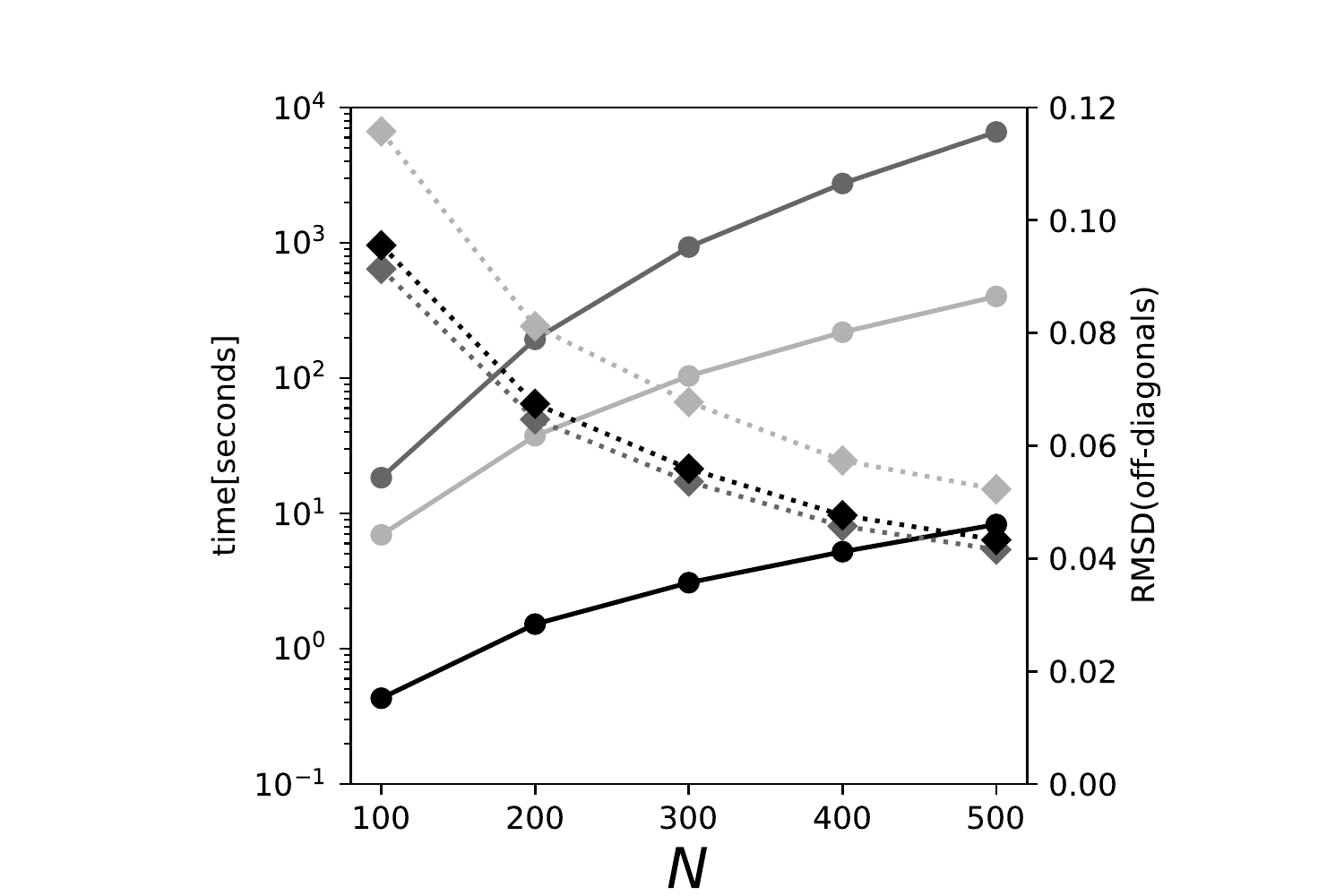}
  &
  \includegraphics[scale=0.4,trim=0.7in 0in 0.7in 0.35in,clip]{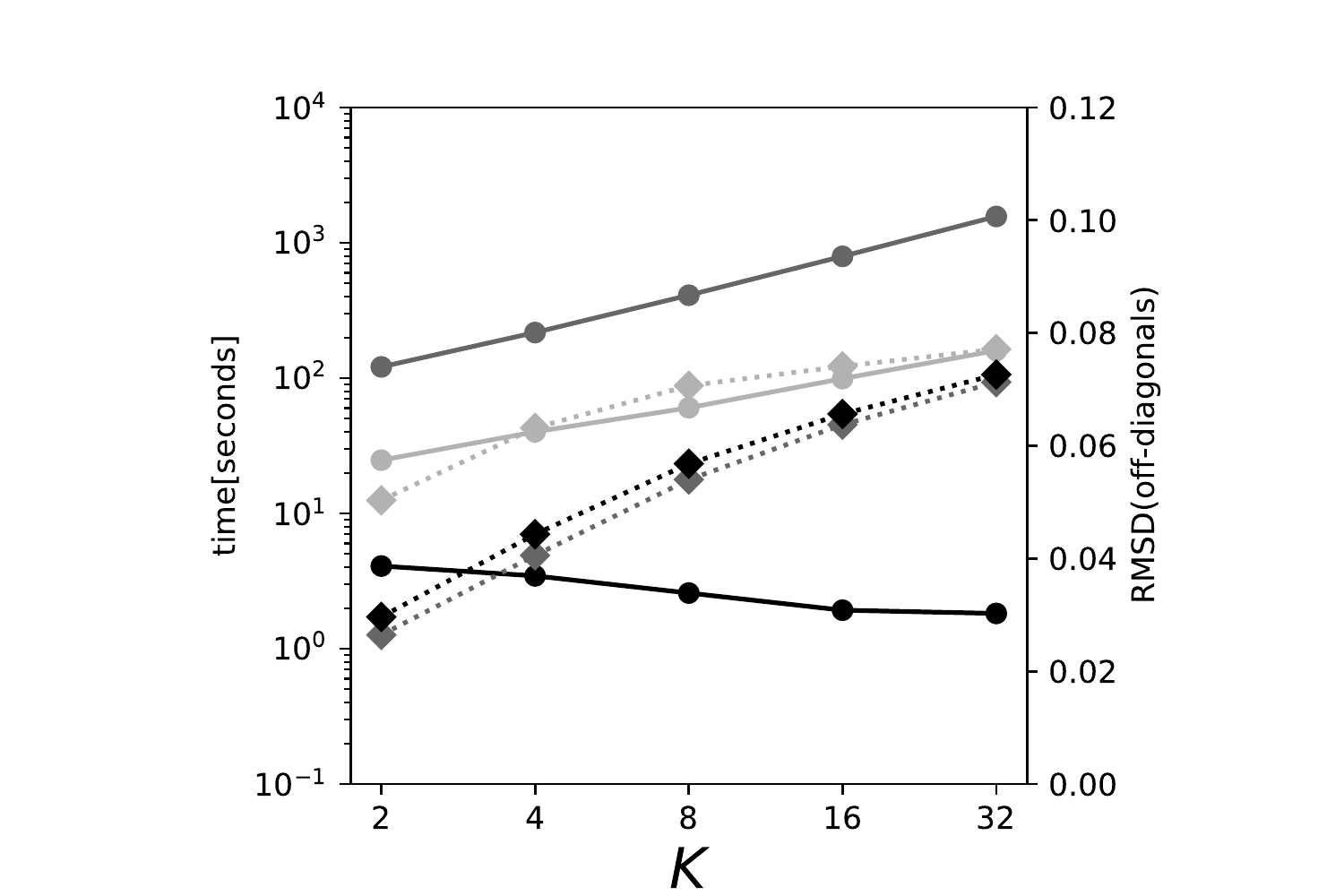} \\ \hline
  \multicolumn{2}{c}{$\alpha=0.25$} \\
  \includegraphics[scale=0.4,trim=0.7in 0in 0.7in 0.35in,clip]{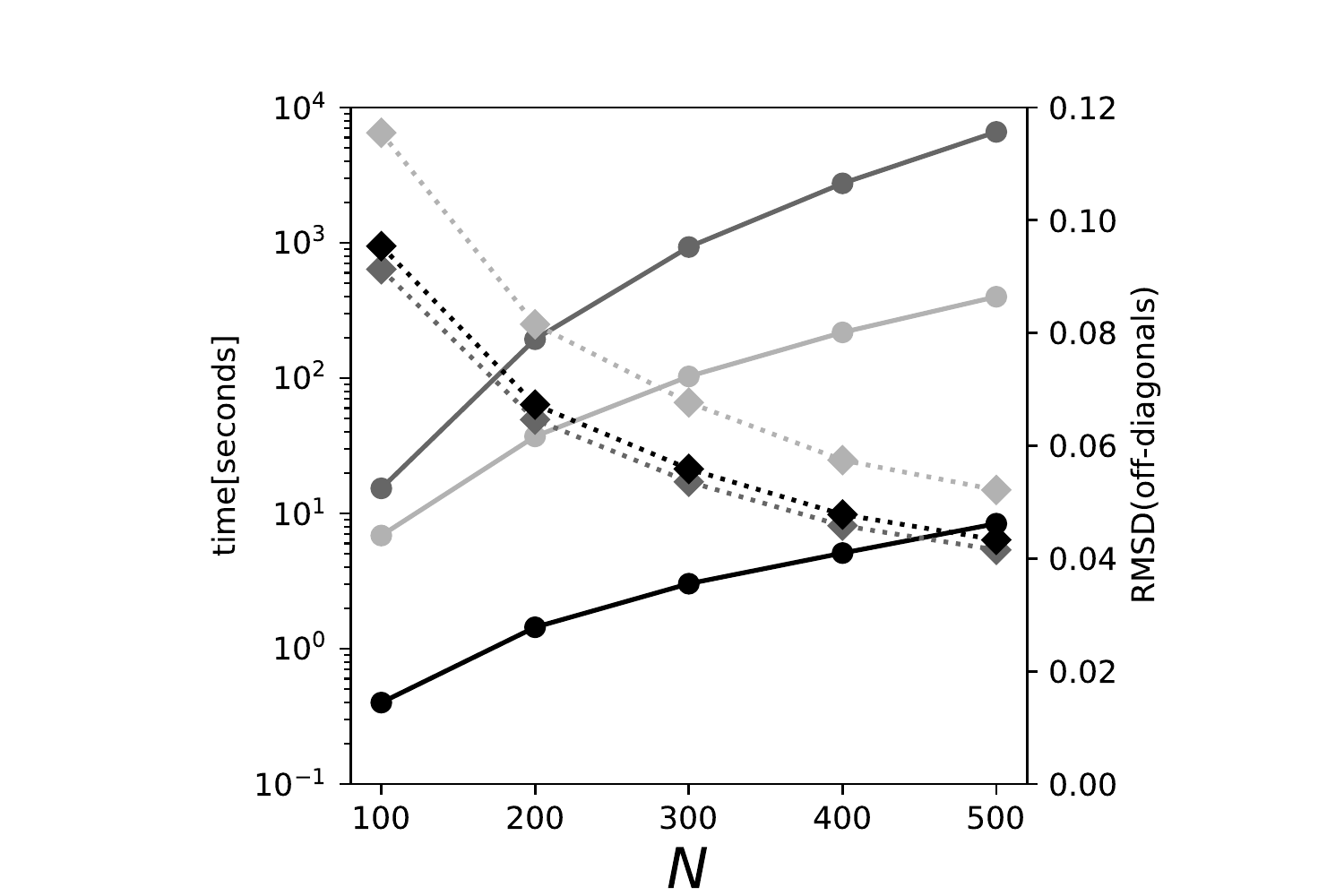}
  &
  \includegraphics[scale=0.4,trim=0.7in 0in 0.7in 0.35in,clip]{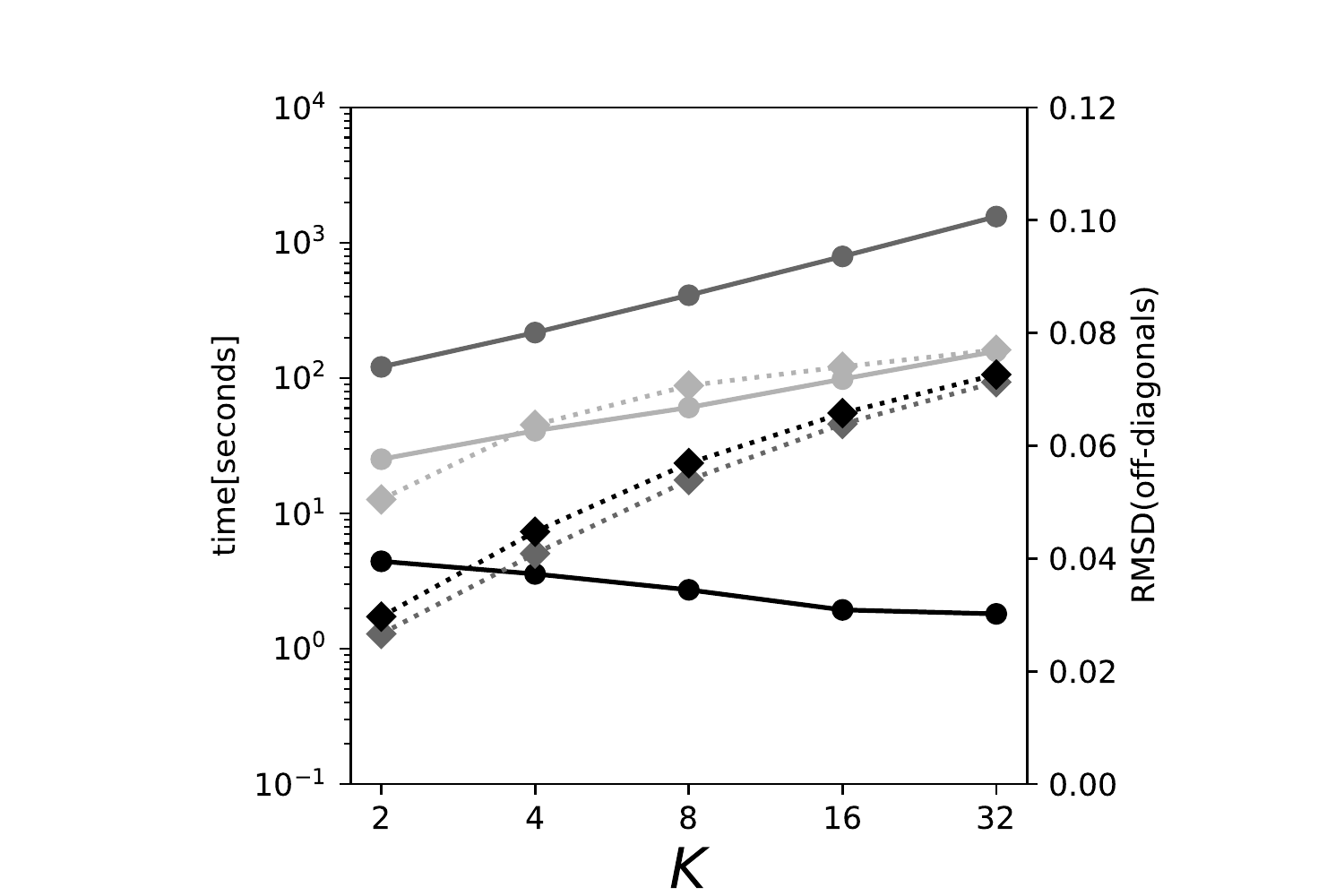} \\ \hline
  \multicolumn{2}{c}{$\alpha=0.50$} \\
  \includegraphics[scale=0.4,trim=0.7in 0in 0.7in 0.35in,clip]{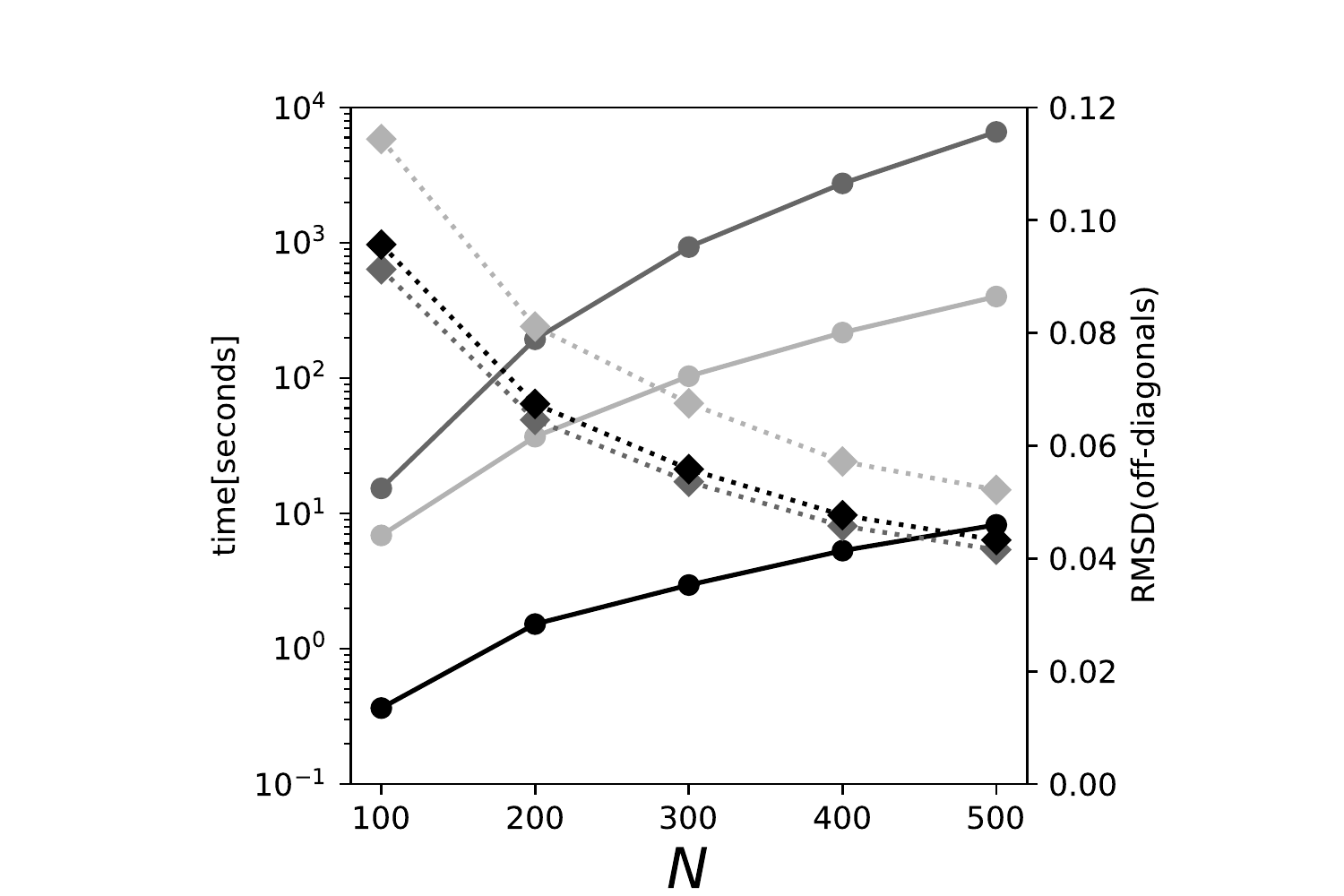}
  &
  \includegraphics[scale=0.4,trim=0.7in 0in 0.7in 0.35in,clip]{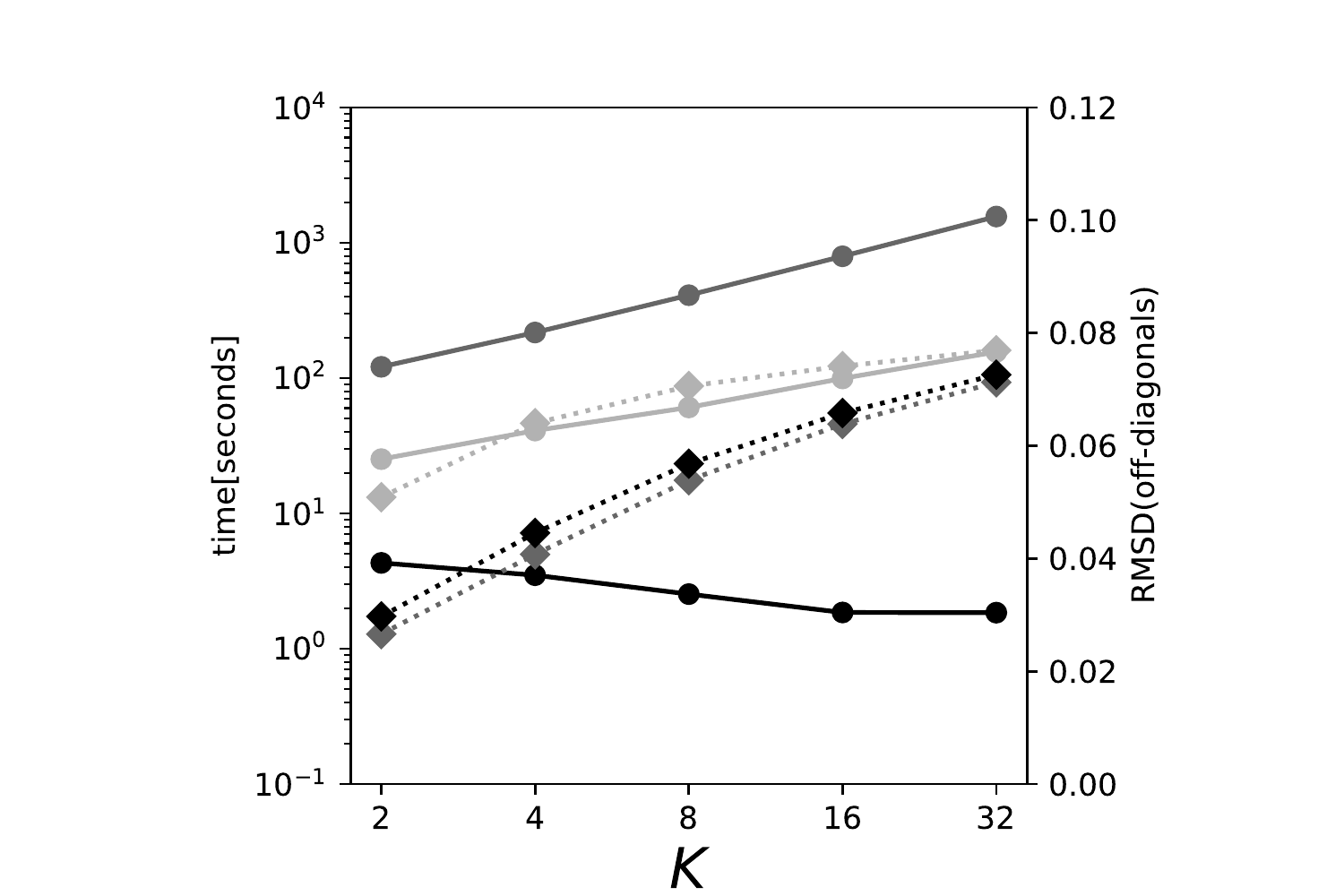} \\ \hline
  \multicolumn{2}{c}{$\alpha=0.75$} \\
  \includegraphics[scale=0.4,trim=0.7in 0in 0.7in 0.35in,clip]{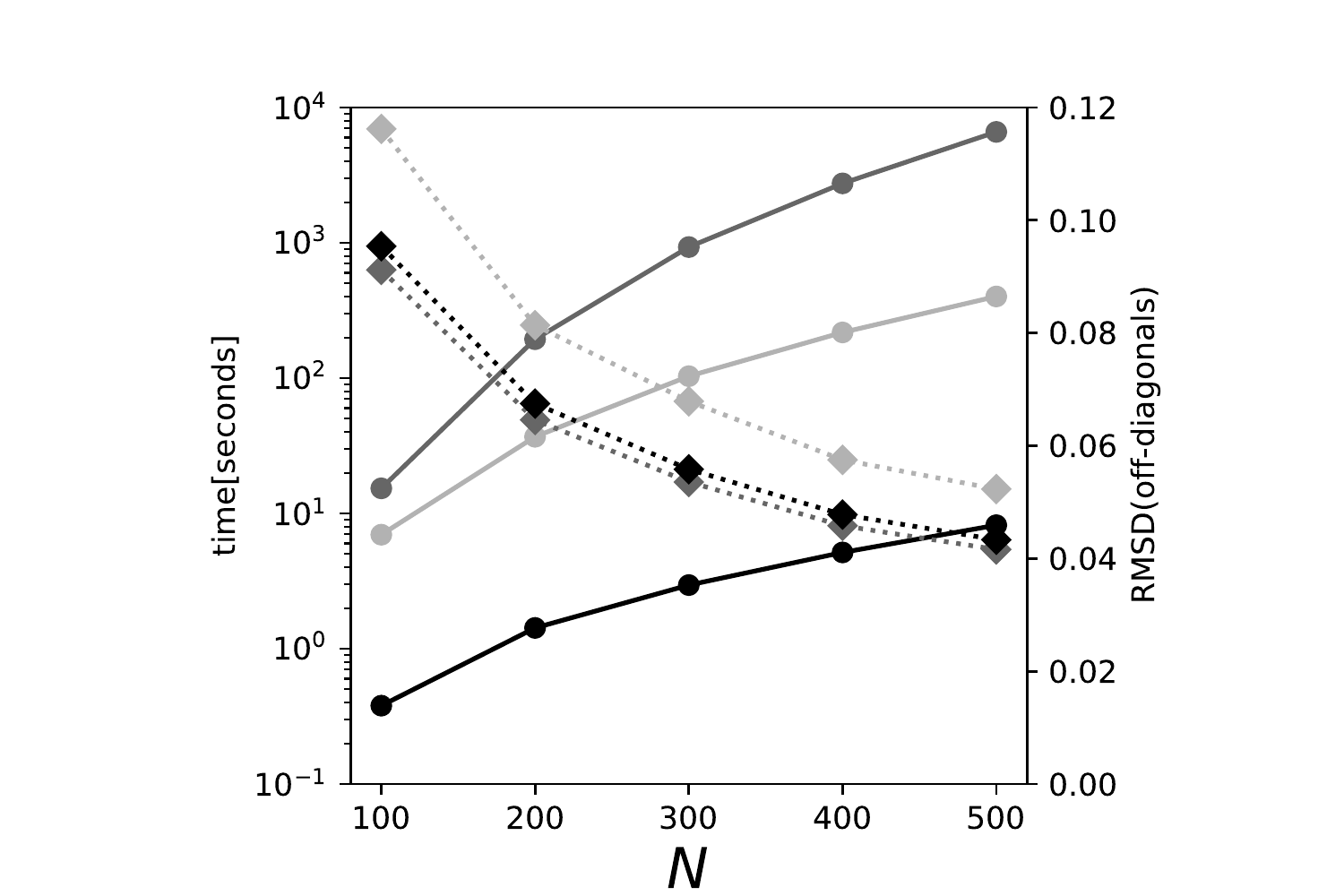}
  &
  \includegraphics[scale=0.4,trim=0.7in 0in 0.7in 0.35in,clip]{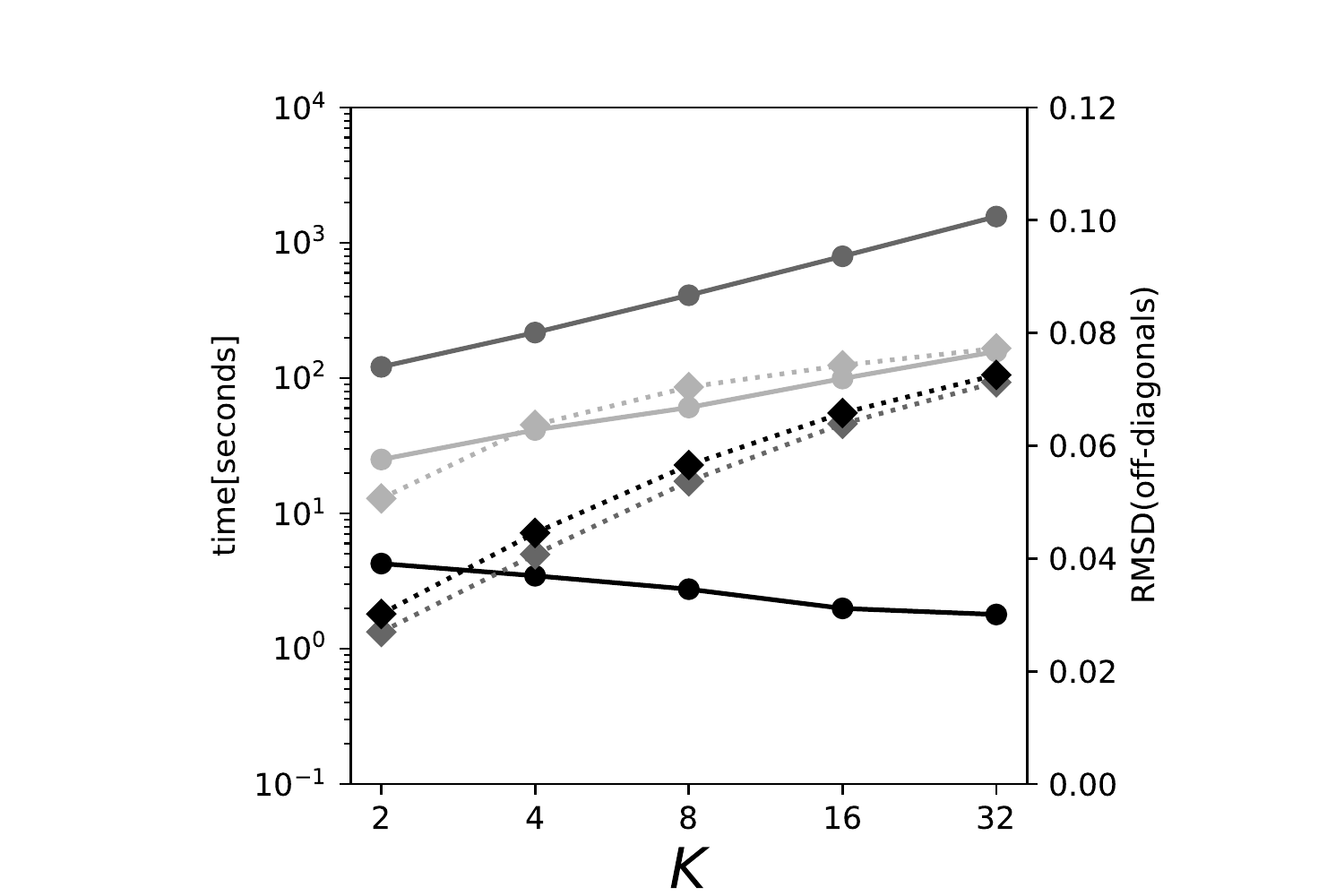} \\ \hline  
  \end{tabular}
  \caption{\label{fig:results}{\scriptsize{Runtime in seconds (solid lines with circles) and root-mean-square deviation (RMSD) from zero of off-diagonal elements (dotted lines with diamonds). Shown is the median across 10 replicates, for various sample sizes $N$ with $K=10$ (left column) and for various number of matrices to diagonalize $K$ with $N=256$ (right column), for various degrees of similarities in eigenvectors across matrices $\alpha$ (rows). Each plot shows lines for JADOC (black), JADE (dark gray), and \texttt{qndiag} (light gray).}}}
\end{figure}

Figure \ref{fig:results} comprises eight line plots (one for each scenario), where each plot shows the median runtime of each algorithm on the left $y$-axis and the median RMSD of each algorithm on the right $y$-axis. Each plot has three lines, corresponding to JADOC, JADE, and \texttt{qndiag}. In the plots on the left (Design 1) $N$ is shown on the $x$-axis and in the plots on the right (Design 2) $K$ is shown on the $x$-axis. In the consecutive rows, results are shown for $\alpha=0,0.25,0.5,0.75$.

We see that JADOC outperforms \texttt{qndiag} by several orders of magnitude in terms of CPU time both as $N$ and as $K$ increase. In turn, \texttt{qndiag} outperforms JADE by another order of magnitude timewise. Overall, JADOC is much fast than \texttt{qndiag} which in turn is much faster than JADE. Also, as expected, for fixed $N$ we see the runtime of JADOC does not increase $K$. In fact, the runtime of JADOC even seems to decrease slightly as $K$ increases.

In terms of the degree of diagonalization, as measured by the RMSD, we see that JADE performs best, very closely followed by JADOC, and finally followed at some distance by \texttt{qndiag}. Thus, in terms of the degree of diagonalization, JADOC closely trails the golden standard provided by JADE, whereas \texttt{qndiag} provides a considerably poorer result.

Overall, the improved runtime afforded by JADOC simply dwarfs the minor increase in the degree of diagonalization that is achieved by JADE.

\section{Conclusions}
\label{sec:conclusions}
We have derived a new framework for finding a square matrix $\x{B}$ such that $\x{B}\x{C}_k\xt{B}$ is as diagonal as possible for $k=1,\ldots,K$, where $\x{C}_1,\ldots,\x{C}_K$ are symmetric, positive (semi)-definite matrices, where $\x{B}\xt{B} = \x{I}_N$ (i.e., $\x{B}$ is orthonormal). Our approach for finding $\x{B}$ we refer to as JADOC: Joint Approximate Diagonalization under Orthogonality Constraints.

In a nutshell, JADOC (\emph{i}) applies a simple dimensionality-reduction technique combined with mild regularization, (\emph{ii}) uses a quasi-Newton approach which is combined with a golden section, to jointly update $\x{B}$ in each iteration, (\emph{iii}) guarantees an orthonormal $\x{B}$, and (\emph{iv}) requires only $O(N^3)$ time per iteration. JADOC is implemented as an open-access tool for Python 3.$x$ and is available on GitHub (see: \url{https://github.com/devlaming/jadoc}).

In most real-life scenarios, positive (semi)-definite matrices $\x{C}_1,\ldots,\x{C}_K$ will not be jointly diagonalizable. In other words, their eigenvectors will not be completely identical. Under this scenario, we find that JADOC improves upon two important approaches, \textit{viz}., the classical JADE method and a recent quasi-Newton approach with a Python package called \texttt{qndiag}. JADOC has far superior runtime compared to both. Moreover, in terms of the degree of diagonalization, JADOC closely follows the golden standard set by JADE.

Finally, we observe, in line with our theory, that for fixed $N$, increasing $K$ does not increase the computational complexity of JADOC. Its time performance remains steady for all values of $K$ considered here. In fact, we even see a slight decrease in CPU time as $K$ increases, as the parallelization implemented in JADOC is tailored towards higher values of $K$.

We conclude that JADOC is an efficient and accurate tool to diagonalize $K$ different $N \times N$, symmetric, positive (semi)-definite matrices, requiring only $O(N^3)$ time per iteration, which constitutes an $O(K)$ reduction in complexity, compared to the $O(KN^3)$ time required in each iteration by competing methods. This reduction in computational complexity makes JADOC far more scalable and, therefore, much more useful than its competitors in the age of big data.

\section*{Acknowledgments}
We would like to acknowledge the fruitful discussions we have had about this topic with Niels Rietveld, Robert Kirkpatrick, and Patrick Groenen. This work was carried out on the Dutch national e-infrastructure with the support of SURF Cooperative (NWO Call for Compute Time EINF-403 to E.A.W.S.).

\bibliographystyle{apalike}
\bibliography{references}

\end{document}